\documentclass[12pt]{article}
\usepackage{geometry}                % See geometry.pdf to learn the layout options. There are lots.
\usepackage{graphicx}
\usepackage{amssymb}
\usepackage{epstopdf}
\newcommand{\vect}[2]{\left(\begin{array}{c}
#1 \\
#2 \\
\end{array}\right)}
\newcommand{\vecd}[3]{\left(\begin{array}{c}
#1 \\
#2 \\
#3 \\
\end{array}\right)}

\title{Affine Maps and Feuerbach's Theorem}
\author{Patrick Morton}
\date{}                                           % Activate to display a given date or no date

\begin{document}
\maketitle
\section{Introduction.}

\noindent In this paper we will study several important affine maps in the geometry of a triangle, and use them to give a surprising proof of Feuerbach's theorem that the incircles and excircles of a triangle are tangent to the nine-point circle.  Recall that the nine-point circle $\textsf{N}$ of a triangle $ABC$ is the circle that lies on the midpoints of the sides, the feet of the altitudes, and the midpoints of the segments joining the orthocenter (the common intersection of the altitudes) to the vertices.  The incircle $\textsf{I}$ is the unique circle which is internally tangent to the three sides of the triangle, while there are three excircles that are externally tangent to all three (extended) sides.  (See Figure 2.)  \medskip

\noindent I will prove that the  that the incircle $\textsf{I}$ and the nine-point circle $\textsf{N}$ are tangent to each other by constructing a mapping which takes $\textsf{N}$ to $\textsf{I}$ and has a fixed point $Z$ lying on both circles.  To show that the circles are tangent to each other at $Z$ our mapping must be a {\it homothety}.  A mapping $h$ on the plane is called a homothety (or dilatation) if there is a point $O$ and a real number $k$ such that: $h(O)=O$; and for any $P \neq O$, the vectors $\overrightarrow{Oh(P)}$ and $\overrightarrow{OP}$ satisfy $\overrightarrow{Oh(P)} = k*\overrightarrow{OP}$, so that $h(P)$ lies on the line $OP$.  A homothety $h$ has the property that for three non-collinear points $O, P, Q$, the triangle $h(OPQ)$ is similar to $OPQ$ with similarity ratio $Oh(P)/OP = |k|$.  Moreover, a line $l$ not through $O$ and its image $h(l)$ under $h$ are always parallel.  The point $O$ is called the center of the homothety. (See [5, p. 105].)  \medskip

\noindent For any triangle $ABC$ there is an important homothety $K$ for which $O=G$, the centroid of $ABC$, and $k=-1/2$.  Since the centroid is two-thirds of the way along the segment joining each vertex to the opposite midpoint, the segment $GA$ has twice the length of $GD_0$, where $D_0$ is the midpoint of side $BC$, and $A$ and $D_0$ are on opposite sides of $G$.  Hence, $K(A)=D_0$.  Thus, $K(ABC) = D_0E_0F_0$, where $E_0$ and $F_0$ are the midpoints of sides $AC$ and $AB$.   The map $K$ is called the complement map.  \medskip 

\noindent The map we will use to prove Feuerbach's Theorem depends on two related affine maps, defined as follows.  An affine map is a map which can be given in terms of Cartesian coordinates as

$$T(x, y) = M {\vect x y} + {\vect {b_1} {b_2}}, $$ \smallskip

\noindent where $M$ is a $2 \times 2$ non-singular matrix with real coefficients and $b_1, b_2 \in {\mathbb R}$.  It is an elementary exercise in linear algebra to show that for any two triangles $ABC$ and $DEF$ there is a unique affine mapping $T$ which maps the vertices of $ABC$ to the vertices of $DEF$.  This result is sometimes referred to as the Fundamental Theorem of Affine Geometry.  (See [2, pp. 53, 69].) Affine maps take lines to lines and preserve parallelism and ratios along parallel lines. The set of affine maps forms an important group called the affine group of the plane.  It is easy to see that any homothety, in particular the map $K$ defined above, is an element of the affine group.  (These facts hold in any Hilbert plane which satisfies the parallel postulate.  In general, $\mathbb{R}$ is replaced by the field of segments of the plane.  See [8].) \medskip

\noindent The first map we need is the unique affine map  $T_1$ taking $ABC$  to $DEF$, where $D, E, F$ are the points of tangency of the incircle $\textsf{I}$ with the respective sides $BC, AC$, and $AB$:

$$T_1(ABC)=DEF.$$ 

\noindent The second map $T_2$ is defined in the following way.  Let $A', B', C'$ be the centers of the three circles which are externally tangent to the three sides of $ABC$, where $A'$ is the center opposite $A$, etc., meaning that the excircle with center $A'$ is tangent to side $BC$ at a point between $B$ and $C$.  The map $T_2$ is then defined by
$$T_2(A'B'C')=ABC.$$ 

\noindent The map $T_2$ could equivalently be defined by $T_2(ABC)=D'E'F'$, where the points $D', E', F'$ are the reflections of the points $D, E, F$ in the midpoints $D_0, E_0, F_0$ of the respective sides $BC, AC, AB$. (See Figure 1 and Lemma 2 in Section 2.)  \medskip

\begin{figure}[htp]
\centering
\includegraphics[height=6in]{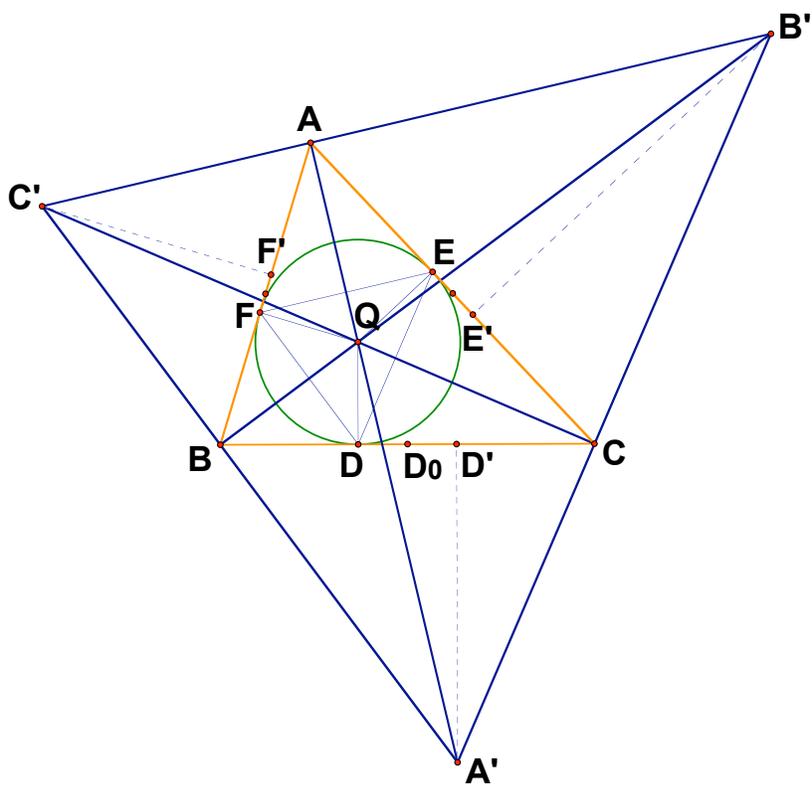}
\caption{Incircle, incenter, and excenters}
\end{figure}

\noindent We are now ready to define the map which will play the main role in our proof.  We set

$$\Phi = T_1K^{-1}T_2K^{-1}, \eqno{(1.1)}$$ \smallskip

\noindent so that $\Phi$ is the composition of the mappings $T_1$ and $T_2$ defined above, alternating with the anti-complement map $K^{-1}$.  The map $\Phi$ is certainly an affine map, but more is true. \bigskip

\noindent {\bf Main Theorem.} 1) The map $\Phi$ takes the nine-point circle $\textsf{N}$ to the incircle $\textsf{I}$. \smallskip

2) $\Phi$ has a fixed point $Z$ which lies on $\textsf{N}$.  \smallskip

3) $\Phi$ is a homothety. \smallskip

4) If the triangle $ABC$ is not equilateral, the unique fixed point (center) of $\Phi$ is the Feuerbach point $Z$, the point where the incircle touches the nine-point circle.   \bigskip

\noindent Once we have proven parts 1)-3) of this theorem, the proof of 4) is easy.  This is because an affine map takes tangent lines to a conic to tangent lines of the image conic.  Since $\Phi$ is a homothety, and $Z$ is fixed, $\Phi$ will fix the tangent line $t$ to $\textsf{N}$ at $Z$.  This implies that $t$ is also the tangent line to $\textsf{I}$ at $Z$, proving that $\textsf{N}$ and $\textsf{I}$ are tangent to each other at $Z$. \medskip

\noindent We will prove an analogous theorem for the excircles of $ABC$.  This yields Feuerbach's famous theorem (see Figure 2): \bigskip

\noindent {\bf Feuerbach's Theorem.}  If the triangle $ABC$ is not equilateral, the nine-point circle $\textsf{N}$ of triangle $ABC$ is tangent to the incircle and to each of the excircles of $ABC$. \bigskip

\noindent When $ABC$ is equilateral it is clear that the nine-point circle coincides with the incircle, and strictly speaking cannot be tangent to itself.  This is the reason for the restricting hypothesis in the above theorems.  \medskip

\noindent The proof we give will show that Feuerbach's theorem holds in any Hilbert plane satisfying the parallel postulate (see Section 6).  In the language of [8], this is equivalent to saying that the theorem holds for the cartesian plane over any ordered, Pythagorean field $\textsf{F}$.  One way to interpret this result is to say that the coordinates of the point of tangency of $\textsf{N}$ and $\textsf{I}$ are in the same field in which the coordinates of the vertices and incenter lie.  This is surprising in itself, since to find the intersection of two circles one generally has to adjoin a square-root.  For other proofs of Feuerbach's Theorem see [1, pp. 105-107] and [5, pp. 136-137].  \medskip

\noindent {\bf Acknowledgements.} I am very grateful to an anonymous referee for his detailed comments on an earlier version of this paper.  His thoughtful analysis gave me additional insight on the proof given here and strongly influenced its presentation.  I am also grateful to Igor Minevich for a careful reading of the manuscript.

\begin{figure}[htp]
\centering
\includegraphics[height=5in]{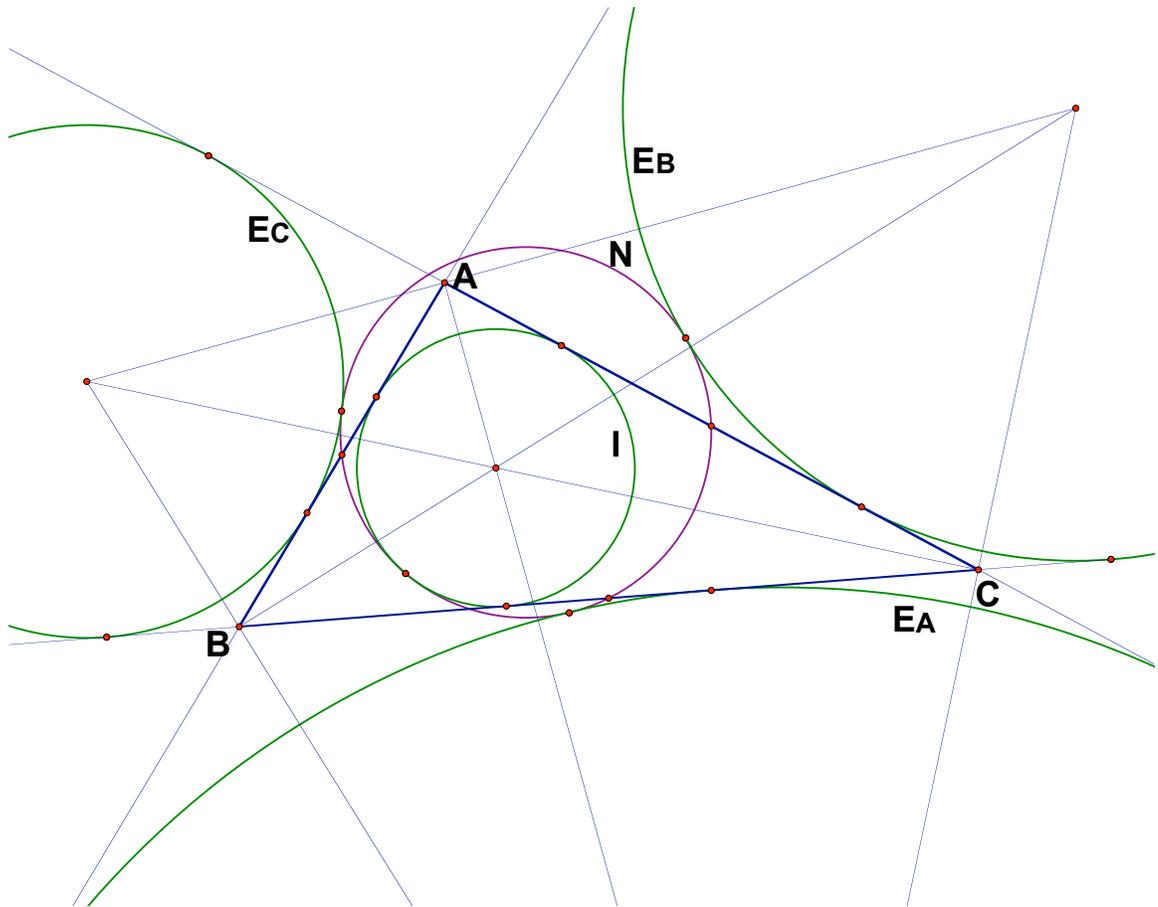}
\caption{Feuerbach's Theorem}
\end{figure}

\section{Proof of Part 1) of the Main Theorem.}

\noindent In this section we give a synthetic proof that $\Phi=T_1K^{-1}T_2K^{-1}$ maps the nine-point circle to the incircle. \medskip

\noindent As in Section 1, we let \textsf{I} denote the incircle of triangle $ABC$, $Q$ its center, and $D, E, F$ the points of tangency of this circle with the sides $BC, AC$, and $AB$, respectively.  Let $T_1$ be the unique affine map taking $ABC$ to $DEF$.  Further, let $A', B', C'$ denote the excenters of $ABC$ opposite the respective vertices $A, B, C$, and let $T_2$ be the affine map taking $A'B'C'$ to $ABC$.  We first consider the map $S=T_1T_2$, which takes triangle $A'B'C'$ to triangle $DEF$.  \medskip

\noindent We note that $Q$ is the orthocenter of triangle $A'B'C'$.  This is because $AQ=A'Q$ is the internal angle bisector of the angle $A$,  and $B'C'$ is the external angle bisector at $A$, so that $A'Q  \bot  B'C'$.  Thus, $A'Q$ is the altitude from vertex $A'$, and in the same way $B'Q=BQ$ and $C'Q=CQ$ are altitudes.  On the other hand $AQ$ is the internal angle bisector in isosceles triangle $EAF$, so that $AQ \bot EF$.  Hence, $EF \| B'C'$.  In the same way, $DE \| A'B'$ and $DF \| A'C'$. \medskip

\noindent Now we appeal to a special case of Desargues' Theorem.  If $A'B'C'$ and $DEF$ are triangles for which $EF \| B'C'$, $DE \| A'B'$, and $DF \| A'C'$, then the triangles are perspective, which means that lines joining corresponding vertices are either all parallel to each other or they are concurrent.  In this case $B'E$ and $C'F$ clearly intersect inside triangle $DEF$, so $A'D$, $B'E$, and $C'F$ must meet at some point $J$.  Now the homothety with center $J$ and ratio $k = \frac{JD}{JA'}$ maps $A'B'C'$ to $DEF$ (by similar triangles) and so must coincide with the map $S$, by the Fundamental Theorem of Affine Geometry.  This proves the following lemma.   \bigskip

\noindent {\bf Lemma 1.} The map $S=T_1T_2$ is a homothety.  \bigskip

\noindent This shows that $T_1$ and $T_2$ are inverse maps on the line at infinity, i.e. on slopes of lines.  (See [5, pp. 56-57] or [6, pp. 81-83].)  Now the map $K^{-1}$ leaves slopes of lines invariant, so $\Phi$ has the same property.  It follows that $Y\Phi(Y)$ is an invariant line, for any ordinary point $Y$ (which is not fixed by $\Phi$).  If two such lines intersect, then $\Phi$ is an affine map that has a fixed point and takes lines to parallel lines, so that $\Phi$ is a homothety.  If all lines of the form $Y\Phi(Y)$ are parallel, then $\Phi$ must be a translation.  Thus, $\Phi$ is either a homothety or a translation. (See [5, p. 120].)  \medskip

\noindent For the proofs in Sections 4 and 5 we need to know the effect of the map $T_2$ on the triangle $ABC$.  \bigskip

\noindent {\bf Lemma 2.} We have $T_2(ABC)=D'E'F'$, where:
 $D'$ is the reflection of $D$ across the midpoint $D_0$ of $BC$, 
 $E'$ is the reflection of $E$ across the midpoint $E_0$ of $AC$,
and $F'$ is the reflection of $F$ across the midpoint $F_0$ of $BC$.  \medskip
   
\noindent {\it Proof.}  We shall prove that $T_2(A)=D'$ by comparing the ratios $C'A/C'B'$ and $CD'/CB$ and using the fact that $T_2(B') = B$ and $T_2(C') = C$.  We first note that the circle $\Gamma$ on $B'$ and $C'$ whose center is the midpoint $M$ of $B'C'$ lies on the points $B$ and $C$, since the angles $\angle C'CB'$ and $\angle B'BC'$ are right angles.  (In other words, $B'C'CB$ is a cyclic quadrilateral.)  Hence, the angles $\angle B'C'Q$ and $\angle CBQ$ which cut off the chord $B'C$ of $\Gamma$ are congruent, as are angles $\angle C'B'Q$ and $\angle BCQ$.  It follows that $B'C'Q \sim CBQ$.  Moreover, $QA$ and $QD$ are altitudes of $B'C'Q$ and $CBQ$, respectively, so that $C'AQ \sim BDQ$.  This gives that

$$\frac{C'A}{BD}=\frac{C'Q}{BQ}=\frac{C'B'}{BC},$$ \smallskip

\noindent and hence that 

$$\frac{C'A}{C'B'}=\frac{BD}{BC}=\frac{CD'}{CB}$$ \smallskip

\noindent since $BD \cong CD'$ by the definition of the point $D'$.  Using the property that affine maps preserve ratios along lines [2] and that $T_2(A)$ lies between points $B$ and $C$, this implies that $T_2(A)=D'$.  In the same way we get that $T_2(B)=E'$ and $T_2(C)=F'$.  \bigskip

\noindent {\bf Remark.}  It can be shown that the point $D'$ is the point of tangency with side $BC$ of the excircle whose center is $A'$, with analogous statements for $E'$ and $F'$.  See [1, p. 160]. \bigskip

\noindent {\bf Theorem 1.} If $\Phi=T_1K^{-1}T_2K^{-1}$, then $\Phi(\textsf{N})=\textsf{I}$. \medskip

\noindent {\it Proof.}  We consider the result of applying the individual maps that define $\Phi$ to the nine-point circle $\textsf{N}$ of $ABC$.  \medskip

\noindent {\it Step 1.} The circle  $\textsf{N}$ passes through the midpoints of the sides of $ABC$, so $K^{-1}( \textsf{N})$ is a circle lying on the vertices of $ABC$.  Hence, $K^{-1}( \textsf{N})=\textsf{C}$ is the circumcircle of $ABC$, whose center is the circumcenter $O$.  The circumcircle $\textsf{C}$ is also the nine-point circle of the excentral triangle $A'B'C'$, because it lies on $A, B$, and $C$, which are the feet of the altitudes of this triangle.  It follows that $\textsf{C}$ passes through the midpoints of the sides of $A'B'C'$.  \medskip

\noindent {\it Step 2.} We apply the map $T_2$, using the fact that the image of a circle under an affine map is an ellipse [2, p. 85].  The ellipse $\textsf{E}_1=T_2K^{-1}( \textsf{N})=T_2(\textsf{C})$ is a circumconic of the triangle $T_2(ABC)=D'E'F'$.  Since affine maps preserve ratios of lengths along a line, the image $T_2(O)$ of the center $O$ is the center of the ellipse.  Moreover, $\textsf{E}_1=T_2(\textsf{C})$ passes through the midpoints $D_0, E_0, F_0$ of the sides of $T_2(A'B'C')=ABC$!  Hence, $\textsf{E}_1$ is an ellipse on the six points $D', E', F', D_0, E_0, F_0$.  (Some of these points may coincide with one another for certain triangles $ABC$.) \medskip

\noindent {\it Step 3.} Similarly, the image $K^{-1}(\textsf{E}_1)=K^{-1}T_2(\textsf{C})$ is an ellipse $\textsf{E}_2$ with center $K^{-1}T_2(O)$.  It passes through the points $K^{-1}(D_0)=A, K^{-1}(E_0)=B, K^{-1}(F_0)=C$, so $\textsf{E}_2$ is a circumconic of $ABC$. \medskip

\noindent {\it Step 4.}  Finally, $\textsf{E}_3=T_1(\textsf{E}_2)$ is an ellipse on the points $T_1(A) = D, T_1(B) = E, T_1(C) = F$.  If we show that  $\textsf{E}_3$ is a circle, then we are done, since the unique circle on the points $D, E, F$ is the incircle  $\textsf{I}$.  But this is easy, because $ \textsf{E}_3 = \Phi( \textsf{N}).$  We have shown that $\Phi$ is either a homothety or a translation.  Therefore, $\textsf{E}_3$, as the image of the circle $\textsf{N}$, must also be a circle.  Hence, $\Phi(\textsf{N})=\textsf{I}$.  \bigskip

\noindent {\bf Corollary.}  The center $N$ of the nine-point circle $\textsf{N}$ is related to the incenter $Q$ by the formula $\Phi(N)=T_1K^{-1}T_2K^{-1}(N)=Q$.  Equivalently, $K^{-1}(N)=O$ is the circumcenter of $ABC$ and $O=T_2^{-1}KT_1^{-1}(Q)$.  \bigskip

\noindent {\bf Remarks.} 1. It can be shown (independent of Feuerbach's theorem)  that unless triangle $ABC$ is equilateral, the radius $r_n$ of the circle $\textsf{N}$ is larger than the radius $r_i$ of the incircle $\textsf{I}$.  This follows from the fact that $2r_n=R$, where $R$ is the radius of the circumscribed circle [6, p. 407]; and Euler's relation $d^2=R(R-2r_i)=4r_n(r_n-r_i)$ for the distance $d$ between the circumcenter and the incenter of $ABC$ [1, p. 85].  Moreover, we have the formulas $\alpha=sr_i$ and $8\alpha r_n=abc$, where $\alpha=\alpha(ABC)$ is the area of $ABC$; $a, b$, and $c$ are the lengths of the sides; and $s = (a+b+c)/2$ is the semi-perimeter.  See [4, pp. 12-13].  From these formulas and Heron's formula for $\alpha$ it follows that 

$$\det(\Phi)=\left(\frac{r_i}{r_n}\right)^2=\left(\frac{8\alpha^2}{abcs} \right)^2 = \left(\frac{(a+b-c)(a-b+c)(-a+b+c)}{abc} \right)^2$$ \smallskip

\noindent is the factor by which the map $\Phi$ changes areas.  In particular, $\Phi$ cannot be a translation unless $ABC$ is equilateral, in which case $\Phi$ is the identity map.  We will see this more simply in Section 5, once we prove that $\Phi$ has a fixed point. \medskip

2. Lemma 2 shows that the maps $T_1$ and $T_2$ can be defined in a symmetrical way.  In [11, pp. 26-28] it is proved that these two maps are conjugate to each other in the affine group, i.e. $T_2=UT_1 U^{-1}$ for some affine map $U$, whenever triangle $ABC$ is not isosceles.  In this case $\det(T_1)=\det(T_2)$, so using the fact that $\det(K)=1/4$, Remark 1 implies that

$$\det(T_1)=\det(T_2)= \frac{(a+b-c)(a-b+c)(-a+b+c)}{4abc}.$$ \smallskip

\noindent This ratio is therefore the ratio of areas $\alpha(DEF)/\alpha(ABC)=\alpha(D'E'F')/\alpha(ABC)$.  We also note the curious fact that the triangles $T_1(D'E'F')=T_1T_2(ABC)$ and $T_2(DEF)=T_2T_1(ABC)$ are always congruent to each other, since the commutator $T_2T_1T_2^{-1}T_1^{-1}$ is a translation taking $T_1(D'E'F')$ to $T_2(DEF)$.  See [11, pp. 29-31].

\section{The excircles of $ABC$.}

\noindent To prove Feuerbach's Theorem we must prove a theorem analogous to Theorem 1 for the excircles of $ABC$.  Consider the excircle $\textsf{E}_A$ whose center is $A'$.  To avoid confusion with the points defined earlier, we label the center of  $\textsf{E}_A$ as $Q_1(=A')$, and label the other three centers $A_1, B_1, C_1$, where $AA_1, BB_1$, and $CC_1$ all lie on $Q_1$.  Point $A_1$ is the previous $Q$, while $B_1=C'$ and $C_1=B'$.  (See Figure 3.) \medskip

\noindent As before, we take $D_1, E_1, F_1$ to be the points of tangency of $\textsf{E}_A$ with the sides of $ABC$, but this time $D_1$ lies between $B$ and $C$ while $E_1$ and $F_1$ lie outside the triangle on the extended sides $AC$ and $AB$.  The map $T_1$ is again defined by $T_1(ABC)=D_1E_1F_1$, and the map $T_2$ by $T_2(A_1B_1C_1) = ABC$.  \medskip

\noindent The proof of Lemma 1 goes through unchanged.  Thus, $\Phi_{A}=T_1K^{-1}T_2K^{-1}$ is a homothety or translation.  In the proof of Lemma 2 we again have that $B_1C_1CB$ is a cyclic quadrilateral, giving the pairs of similar triangles $B_1C_1Q_1 \sim CBQ_1$ and $C_1AQ_1 \sim BD_1Q_1$ and leading to $T_2(A)=D_1'$, where $D_1'$ is the reflection of $D_1$ in the point $D_0$.  \medskip

\noindent To prove $T_2(B)=E_1'$ use the fact that $A_1AC_1C$ is a cyclic quadrilateral, implying that $\angle Q_1AC \cong \angle Q_1C_1A_1$  and $A_1C_1Q_1 \sim CAQ_1$.  Corresponding altitudes in these obtuse triangles are $Q_1B$ and $Q_1E_1$, so $Q_1BC_1 \sim Q_1E_1A$.  This implies that

$$\frac{C_1B}{AE_1}=\frac{Q_1C_1}{Q_1A}=\frac{C_1A_1}{AC},$$ \smallskip

\begin{figure}[htp]
\centering
\includegraphics[height=6in]{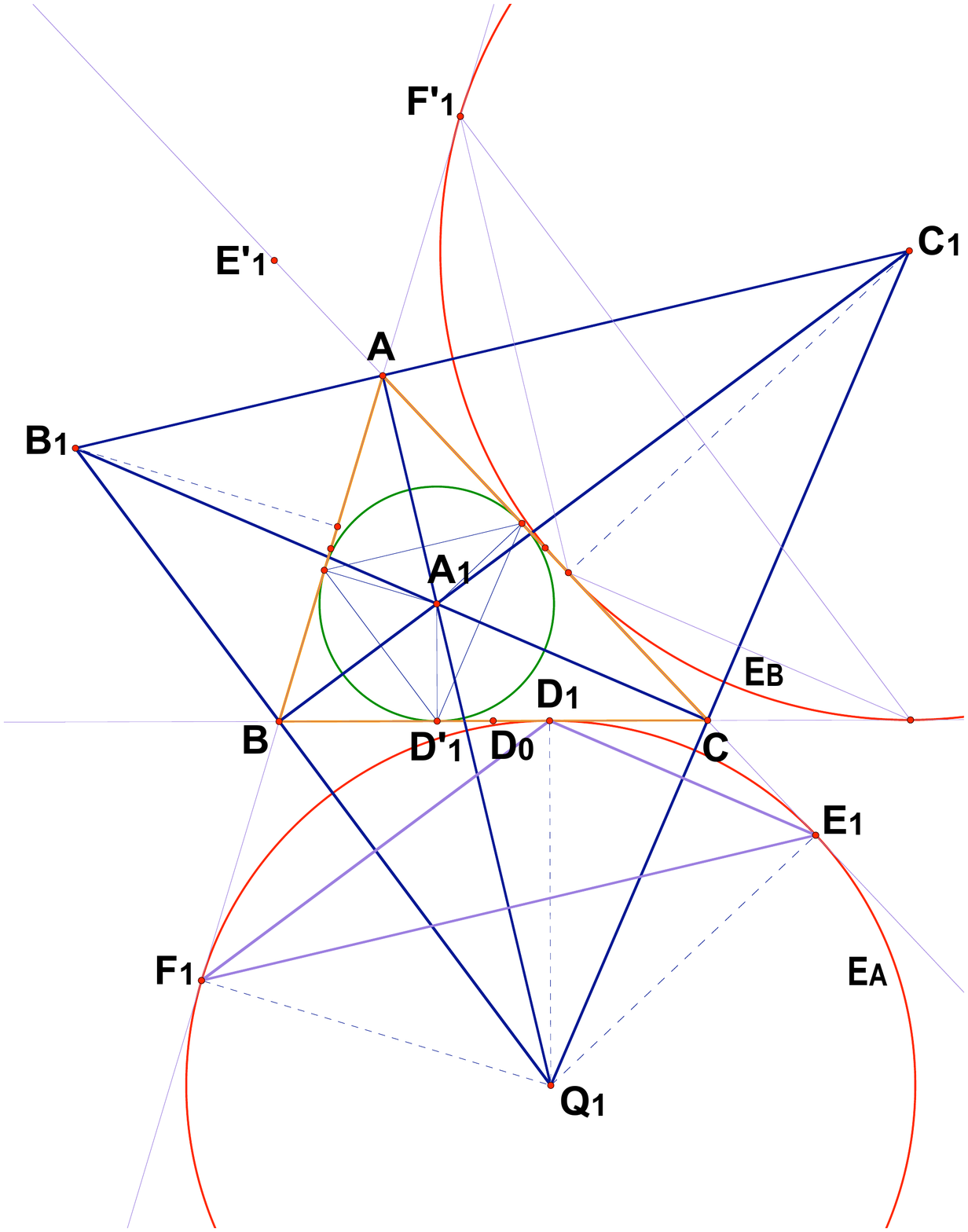}
\caption{Proof that $\Phi_A(\textsf{N})=\textsf{E}_A$}
\end{figure}

\noindent which gives that

$$\frac{C_1B}{C_1A_1}=\frac{AE_1}{AC}=\frac{CE_1'}{CA}.$$ \smallskip

\noindent Now $T_2(A_1)=A$ and $T_2(C_1)=C$.  Since  $A_1$ is between points $B$ and $C_1$, while $A$ is between $E_1'$ and $C$, we conclude that $T_2(B)=E_1'$.  In the same way we get $T_2(C)=F_1'$ and therefore $T_2(ABC)=D_1'E_1'F_1'$, which is the analogue of Lemma 2.  \medskip

\noindent With these additions the proof of Theorem 1 goes through exactly as before, without additional changes. \bigskip

\noindent {\bf Theorem 2.} Define the affine maps $T_1$ and $T_2$ by

$$T_1(ABC)=D_1E_1F_1, \quad T_2(ABC)=D_1'E_1'F_1',$$ \smallskip

\noindent where $D_1, E_1, F_1$ are the points where the excircle $\textsf{E}_A$ touches the respective sides $BC, AC, AB$, and $D_1', E_1', F_1'$ are the reflections of $D_1, E_1, F_1$ in the midpoints $D_0, E_0$, $F_0$.  Then the map $\Phi_{A}=T_1K^{-1}T_2K^{-1}$ takes the nine-point circle $\textsf{N}$ to the excircle $\textsf{E}_A$ with center $Q_1=A'$.  There are corresponding statements for the other two excircles.

\section{Barycentric coordinates and fixed points.}

\noindent We shall now exhibit fixed points of the maps $\Phi$ and $\Phi_A$ by expressing them in terms of homogeneous barycentric coordinates. These two maps are actually two members of a family of maps which can be parametrized in terms of a variable point $P$.  \medskip

\noindent The maps $T_1$ and $T_2$ are determined by the points $D, E, F$ (or $D_1, E_1, F_1$) and their reflections $D', E', F'$.  Since $D, E, F$ are points of tangency of the circle $\textsf{I}$ with the sides of $ABC$, we have $AE \cong AF, BF \cong BD, CD \cong CE$, which implies that

$$\frac{BD}{DC} \cdot \frac{CE}{EA} \cdot \frac{AF}{FB} = 1.$$  \smallskip

\noindent Ceva's Theorem [1, pp. 158-160] implies that the lines $AD, BE$, and $CF$ are concurrent at a point $P$ called the Gergonne point of triangle $ABC$.  A similar argument using  the congruences $BD' \cong CD, CE' \cong AE$, etc., implies that $AD', BE'$, and $CF'$ are also concurrent at a point $P'$.  Thus, the maps $T_1$ and $T_2$ are two instances of the same construction.  In the same way $AD_1, BE_1$, and $CF_1$ are concurrent at an external Gergonne point $P_A$; and $AD_1', BE_1'$, and $CF_1'$ are concurrent at a corresponding point $P_A'$.  \medskip

\noindent We now change our point of view.  We let $P$ be any point not on the sides of $ABC$ or its {\it anticomplementary triangle} $K^{-1}(ABC)$, whose sides are lines through $A, B, C$ parallel to the opposite sides $BC, AC, AB$.  Then $AP$ will intersect side $BC$ in a point $D$; similarly $BP \cdot AC=E$ and $CP \cdot AB=F$, where the dot indicates the intersection of lines.  We again let $D', E', F'$ be the reflections of $D, E, F$ in the points $D_0, E_0, F_0$. For the general point $P$ we define $T_1$, $T_2$, and $\Phi_P$ by

$$T_1(ABC)=DEF, \quad T_2(ABC)=D'E'F', \quad \Phi_P=T_1K^{-1}T_2K^{-1}.$$ \smallskip

\noindent The map $\Phi_P$ coincides with $\Phi$ or $\Phi_A$ when $P$ is taken to be the Gergonne point or the external Gergonne point $P_A$ of triangle $ABC$.  \smallskip 

\noindent We will use a representation of the maps $T_1$ and $T_2$ which arises from {\it homogeneous barycentric coordinates}, as follows.  First, the point $P$ is assigned coordinates

$$P \rightarrow (x,y,z) \hspace{.1 in} \textrm{if} \hspace{.1 in} v_P=\frac{x}{s}v_A+\frac{y}{s}v_B+\frac{z}{s}v_C, \hspace{.1 in} s = x+y+z\eqno{(4.1)}$$ \smallskip

\noindent where $v_P$ is a vector from the origin to the point $P$ and the vectors $v_A,v_B,v_C$ correspond to the points $A,B,C$.   Any ordinary point $P$ has such a representation with $s=1$, because, for example $v_D=\lambda v_B+(1-\lambda)v_C$ for some $\lambda$ and $v_P=\mu v_A+(1-\mu)v_D$ for another constant $\mu$.  Putting these together yields the expression

$$v_P=\mu v_A+\lambda(1-\mu)v_B+(1-\mu)(1-\lambda)v_C,$$ \smallskip

\noindent which has the form (4.1).  If we set $v_P=(u,v)$ and $v_A=(a_1, a_2), v_B=(b_1, b_2), v_C=(c_1, c_2)$, then (4.1) can be written in matrix form as

$${\vecd u v  1} = \left(
\begin{array}{ccc}
{a_1} & {b_1}  & {c_1}  \\
 {a_2} & {b_2}  & {c_2}  \\
 1 & 1 & 1\\
\end{array}
\right)
{\vecd {x/s} {y/s} {z/s}}=\textsf{A}{\vecd {x/s} {y/s} {z/s}}.\eqno{(4.2)}$$ \smallskip

\noindent  The point $P$ is represented by all the coordinate triples $(\lambda x, \lambda y, \lambda z)$ for non-zero $\lambda$ and only by these.  This is because the non-collinearity of the points $A, B$, and $C$ implies that the matrix $\textsf{A}$ in (4.2) is non-singular, so the coordinate triple $(x/s,y/s,z/s)$ is uniquely determined.  Conversely, any triple $(x,y,z)$ with $s=x+y+z \neq 0$ corresponds to a unique ordinary  $P$ in the plane, by the assignment (4.1).   Note that the vertices $A, B, C$ and the centroid $G$ are assigned coordinates

$$A \rightarrow (1, 0, 0), \quad B \rightarrow (0, 1, 0), \quad C \rightarrow (0, 0, 1), \quad G \rightarrow (1, 1, 1).$$ \smallskip

\noindent It is convenient at this point to ``projectivize" equation (4.2) by multiplying through by $\lambda=s$ and writing it in the form

$$\lambda{\vecd u v  1} = \textsf{A} {\vecd {x} {y} {z}}, \quad \lambda \neq 0.\eqno{(4.3)}$$

\noindent In this form we can multiply $x, y, z$ by a common non-zero factor and the equation remains true.  \medskip

\noindent Now what happens if $(x,y,z) \neq (0,0,0)$ is a triple for which $s=x+y+z=0$?  In this case (4.3) becomes

$$\lambda{\vecd u v  0} = \textsf{A} {\vecd {x} {y} {z}}, \quad \lambda \neq 0, \eqno{(4.4)}$$ \smallskip

\noindent where not both $u$ and $v$ are $0$.  Hence, all the vectors
$$\mu (xv_A+yv_B+zv_C)= \mu (xv_A+yv_B-(x+y)v_C)= \mu x(v_A-v_C)+\mu y(v_B-v_C)$$

\noindent lie in the same direction, on a line with slope $m=v/u$.  Such a triple $(x, y, z)$ therefore defines an equivalence class of parallel lines, which we identify with an ``ideal" point on the line at infinity $l_\infty$.  Any class of parallel lines corresponds to a set of vectors of the above form, since $v_A-v_C$ and $v_B-v_C$ are linearly independent.  Hence, the points on $l_\infty$ are the points for which $s=x+y+z=0$ in this coordinate system. (We refer the reader to [4, pp. 216-221] and [12] for more on homogeneous barycentric coordinates and their many uses.)  \medskip

\noindent If $\alpha u + \beta v + \gamma =0$ is the equation of a line $l$ in cartesian coordinates, then by (4.3) the equation of $l$ in barycentric coordinates is

$$(\alpha \quad \beta \quad \gamma) \textsf{A} {\vecd {x} {y} {z}} =0.$$ \smallskip

\noindent Thus lines in the $uv$ plane are represented by homogeneous linear equations in $x, y, z$. In particular, it is easy to see that the equations of the sides of $ABC$ in this coordinate system are simply

$$AB: z=0; \hspace{.2 in} AC: y=0; \hspace{.2 in} BC: x=0.$$ \smallskip

\noindent The key point for us is what affine transformations look like in homogeneous coordinates.  In the affine transformation

$$T{\vect u v} = \left(
\begin{array}{cc}
{a} & {b}  \\
 {c} & {d}  \\
\end{array}
\right) {\vect u v} + {\vect {u_1} {v_1}}$$

\noindent we augment the matrices on the right-hand side and view $T$ as a transformation on points in the plane $w=1$ (in $uvw$ space):

$$T{\vecd u v 1} = \left(
\begin{array}{ccc}
{a} & {b} & 0 \\
 {c} & {d} & 0 \\
 0 & 0 & 1\\
\end{array}
\right) {\vecd u v 1} + {\vecd {u_1} {v_1} 0}=
\left(
\begin{array}{ccc}
{a} & {b}  & {u_1}  \\
 {c} & {d}  & {v_1}  \\
 0 & 0 & 1\\
\end{array}
\right) {\vecd u v 1}.$$ \smallskip

\noindent Now (4.3) gives that

$$\lambda T{\vecd u v 1} = \left(
\begin{array}{ccc}
{a} & {b}  & {u_1}  \\
 {c} & {d}  & {v_1}  \\
 0 & 0 & 1\\
\end{array}
\right) \textsf{A} {\vecd {x} {y} {z}},$$ \smallskip

\noindent so using (4.3) once more shows that $T$, as a mapping on homogeneous barycentric coordinates, is given by a non-singular $3 \times 3$ matrix:

$$ {\vecd {x'} {y'} {z'}} = T{\vecd x y z} = \lambda_1 \textsf{A}^{-1} \left(
\begin{array}{ccc}
{a} & {b}  & {u_1}  \\
 {c} & {d}  & {v_1}  \\
 0 & 0 & 1\\
\end{array}
\right) \textsf{A} {\vecd {x} {y} {z}}, \quad \lambda_1 \neq 0.\eqno{(4.5)}$$ \smallskip

\noindent Furthermore, it is not hard to check that a mapping of the form (4.5) maps the line at infinity $x+y+z=0$ to itself.  Conversely, any non-singular $3 \times 3$ matrix that maps $l_\infty$ to itself represents an affine mapping.  (See [6, pp. 462-463].)  \medskip

\noindent The reason for using homogeneous barycentric coordinates is that calculations become much more manageable, since non-zero homogeneous factors can be ignored or introduced when convenient.  In particular, we can work with polynomial expressions instead of rational expressions.  For example, (4.1) gives that $sv_P-xv_A=yv_B+zv_C$, which implies that the intersection $D=AP \cdot BC$ has coordinates $(0,y,z)$.  Thus we can say that

$$ D \rightarrow (0,y,z), \hspace{.1 in} E \rightarrow (x,0,z), \hspace{.1 in} F \rightarrow (x,y,0).\eqno{(4.6)}$$ \smallskip

\noindent  We use the fact that $BC$ and $DD'$ have the same midpoint to deduce $v_{D'}=v_B+v_C-v_D$, and therefore

$$v_{D'}=v_B+v_C-\bigg(\frac{y}{y+z}v_B+\frac{z}{y+z}v_C\bigg) = \frac{z}{y+z}v_B + \frac{y}{y+z}v_C.$$ \smallskip

\noindent Thus we have

$$D' \rightarrow (0,z,y) = (0,1/y,1/z), \hspace{.05 in} E' \rightarrow (z,0,x)=(1/x,0,1/z),$$

$$ \hspace{.05 in} F' \rightarrow (y,x,0)=(1/x,1/y,0).$$ \smallskip

\noindent From these representations it is easy to see that the point

$$P' \rightarrow (1/x,1/y,1/z) = (yz,xz,xy), \hspace{.2 in} xyz \neq 0,\eqno{(4.7)}$$ \smallskip

\noindent is the intersection of the lines $AD', BE'$ and $CF'$.  This point $P'$ is called the isotomic conjugate of $P$.  \medskip

\noindent The complement map has a very nice representation in terms of barycentric coordinates.  Since $K$ takes the vertices of $ABC$ to the midpoints of the opposite sides, we have for example that $K(A)=D_0=(0,1,1)$.  Thus, we may represent $K$ in homogeneous barycentric coordinates by the matrix

$$K=
\left(
\begin{array}{ccc}
 0 & 1  & 1  \\
 1 & 0  & 1  \\
 1 & 1  &  0 
\end{array}
\right).
\eqno{(4.8)}$$
\smallskip

\noindent It follows easily that the equations of the sides of the anticomplementary triangle $K^{-1}(ABC)$ are given by

$$K^{-1}(BC): y+z=0, \quad K^{-1}(AC): x+z=0, \quad K^{-1}(AB): x+y=0.$$ \smallskip

\noindent In the proofs that follow, we compute the images of various points by identifying them with their homogeneous coordinates and using ordinary matrix arithmetic.  \bigskip

\noindent {\bf Lemma 3.} If $P=(x,y,z)$ does not lie on the extended sides of $ABC$ or its anticomplementary triangle, the matrix representations of the maps $T_1$ and $T_2$ in homogeneous barycentric coordinates are

$$T_1=
\left(
\begin{array}{ccc}
  0 & x'(x+y)  & x'(x+z)  \\
  y'(x+y) & 0  & y'(y+z)  \\
  z'(x+z) & z'(y+z)  & 0  
\end{array}
\right)\eqno{(4.9)}$$

\noindent and
$$T_2 =
\left(
\begin{array}{ccc}
 0 & z'(y+z)  & y'(y+z)  \\
  z'(x+z) & 0  & x'(x+z)  \\
  y'(x+y) & x'(x+y)  & 0  
\end{array}
\right), \eqno{(4.10)}$$
\smallskip

\noindent where $x'=x(y+z), y'=y(z+x), z'=z(x+y)$.  \medskip

\noindent {\it Proof.} We check easily that the map $T_1$ defined by the matrix in (4.9) satisfies

$$T_1(A) = T_1{\vecd 1 0 0} = {\vecd 0 {y(x+y)(x+z)} {z(x+y)(x+z)}} = {\vecd 0 y z} = D,$$ \smallskip

\noindent since $(x+y)(x+z) \neq 0$ for points $P$ not on the sides of the anticomplementary triangle of $ABC$.  Similarly $T_1(B) = E$ and $T_1(C)=F$.  In the same way we see that $T_2(ABC)=D'E'F'$ since $xyz \neq 0$ for points $P$ not on the sides of $ABC$.   It remains only to check that the matrices in (4.9) and (4.10) define affine maps.  For $T_1$ we check this as follows: if $(u,v,w)$ is a point on the line at infinity $l_\infty$, then $u+v+w = 0$ and the sum of the coordinates of the image of this point under $T_1$ is

$$(1,1,1) \cdot T_1{\vecd u v w} = (x+y)(x+z)(y+z) (1,1,1) \cdot {\vecd u v w} $$

$$= (x+y)(x+z)(y+z)(u+v+w) = 0.$$ \smallskip

\noindent Hence, $T_1$ maps $l_\infty$ to itself and is therefore an affine map.  A similar computation shows the same for $T_2$.  \bigskip

\noindent {\bf Remark.}  The matrices $T_1$ and $T_2$ in (4.9) and (4.10) satisfy $T_1+T_2=(x+y)(x+z)(y+z)K$, where $K$ is the complement map.  This means that for any point $X$ the midpoint of $T_1(X)$ and $T_2(X)$ is the point $K(X)$.  Therefore, as affine maps, $T_1$ and $T_2$ satisfy $\frac{1}{2}(T_1+T_2)=K$.  \bigskip

\noindent {\bf Lemma 4.} For any point $P=(x,y,z)$ not on the extended sides of triangle $ABC$ or its anticomplementry triangle, the map $\Phi_P$ fixes the point $Z = (x(y-z)^2,y(x-z)^2,z(x-y)^2)$. \medskip

\noindent {\it Proof.} Using the fact that 

$$K^{-1} = \left(
\begin{array}{ccc}
 -1 & 1  & 1  \\
  1 & -1  & 1  \\
  1 & 1  & -1 
\end{array}
\right),
$$ \smallskip

\noindent we check that $Z = (x(y-z)^2,y(x-z)^2,z(x-y)^2)$ is a fixed point of both maps $T_1K^{-1}$ and $T_2K^{-1}$.  This follows from

$$K^{-1}(Z) = K^{-1}{\vecd {x(y-z)^2} {y(x-z)^2} {z(x-y)^2} } = {\vecd {(y+z)(x-z)(x-y)} {(x+z)(y-x)(y-z)} {(x+y)(z-x)(z-y)}}, $$ \smallskip

\noindent so that the first coordinate of $T_1K^{-1}(Z)$ is

$$x'(x+y)(x+z)(y-x)(y-z)+x'(x+z)(x+y)(z-x)(z-y)$$

$$=(y+z)(x+y)(x+z) \cdot x(y-z)^2.$$ \smallskip

\noindent In the same way we compute the other coordinates and find that $T_1K^{-1}(Z) = Z$.  Now if we substitute the homogeneous coordinates of the point $P'$ for the coordinates of $P$ in the definition of $Z$ we get the same point back again.  Hence $Z$ is symmetrically defined with respect to the coordinates of $P$ and $P'$.  It follows immediately that $Z$ is also a fixed point of $T_2K^{-1}$ and therefore of $\Phi_P$.  Alternatively, we may use the remark after Lemma 3 with $X=K^{-1}(Z)$ to see that the midpoint of $T_1K^{-1}(Z)=Z$ and $T_2K^{-1}(Z)$ is $Z$.  Hence $T_2K^{-1}(Z)=Z$.  This completes the proof.  \bigskip

\noindent {\bf Remark.}  If $P$ is different from the centroid of $ABC$, then the coordinates of $Z$ are not all zero, i.e. $Z$ is a well-defined point in the extended plane.  As an interior point of triangle $ABC$, the coordinates of the Gergonne point $P$ can be taken to be positive, which implies that the coordinates of the corresponding point $Z$ can also be taken to be positive.  In this case $Z$ is an ordinary point of the plane, i.e., it does not lie on $l_\infty$.  It is easy to see that the Gergonne point coincides with the centroid if and only if triangle $ABC$ is equilateral.

\section {Proof of the Main Theorem.}

\noindent To prove part 2) of the Main Theorem we compute an equation for the ellipse $\textsf{E}_1$ in Step 2 of the proof of Theorem 1, given that it passes through the points $D_0, E_0, F_0, D', E', F'$.  We may do this using undetermined coefficients, and we find the equation

$$\textsf{E}_1: F(X,Y,Z) = xX^2 + yY^2 + zZ^2 -(x + y)XY -(x + z)XZ -(y + z)YZ = 0. $$ \smallskip

\noindent For example, the point $D_0 = (0, 1, 1)$ clearly lies on $\textsf{E}_1$, as do $E_0$ and $F_0$; and $D' = (0, z, y)$ lies on this conic because

$$F(0,z,y)=yz^2+zy^2-(y+z)zy = 0.$$ \smallskip

\noindent We show now that the point $Z=(x(y-z)^2,y(x-z)^2,z(x-y)^2)$ lies on $\textsf{E}_1$.  This is a straightforward algebraic calculation, which we arrange as follows.  We first let

$$R=x\cdot x^2(y-z)^4+y \cdot y^2(x-z)^4+z \cdot z^2(x-y)^4,$$ \smallskip
\noindent and
$$S=(x+y) \cdot xy(y-z)^2(x-z)^2+(x+z) \cdot xz(y-z)^2(x-y)^2+(y+z) \cdot yz(x-z)^2(x-y)^2,$$ \smallskip

\noindent so that $F(x(y-z)^2,y(x-z)^2,z(x-y)^2) = R-S$.  The binomial theorem gives that

$$R = (y^3+z^3)x^4+(-8y^3z+y^4+z^4+6y^2z^2-8yz^3)x^3$$

$$+(6y^3z^2+6z^3y^2)x^2-(8y^3z^3)x+y^3z^4+z^3y^4,$$ \smallskip

\noindent while the sum of the first two terms in $S$ is

$$(x+y) \cdot xy(y-z)^2(x-z)^2+(x+z) \cdot xz(y-z)^2(x-y)^2$$

$$=(y-z)^2 \big\{(y+z)x^4+(y^2-4yz+z^2)x^3+(-y^2z-yz^2)x^2+2xy^2z^2\big\}.$$ \smallskip

\noindent Adding the last expression to

$$(y+z) \cdot yz(x-z)^2(x-y)^2=$$

$$(y+z) \big\{yzx^4+(-2y^2z-2yz^2)x^3+(y^3z+4y^2z^2+yz^3)x^2+(-2y^3z^2-2z^3y^2)x+y^3z^3\big\}$$ \smallskip

\noindent gives that

$$S = (y^3+z^3)x^4+(-8y^3z+y^4+z^4+6y^2z^2-8yz^3)x^3$$

$$+(6y^3z^2+6z^3y^2)x^2-(8y^3z^3)x+y^3z^4+z^3y^4,$$ \smallskip

\noindent and hence that $R = S$.  Thus $F(Z)=0$ and the point $Z$ lies on the conic $\textsf{E}_1$, as claimed.  \bigskip

\noindent The conic $\textsf{E}_1$ is called the nine-point conic of the {\it quadrangle} $ABCP'$, which is the figure consisting of the six lines joining different pairs of vertices in $\{A, B, C, P'\}$.  See [3, p. 90]. 

\bigskip

\noindent {\it Completion of the Proof of the Main Theorem.} \bigskip

\noindent We first let $P (\neq $ centroid) be the Gergonne point of triangle $ABC$, so that $\Phi_P=\Phi$.  From Step 2 in the proof of Theorem 1 above we have that $\textsf{E}_1=T_2(\textsf{C})$.  Hence, $T_2^{-1}(Z)$ lies on the circumcircle $\textsf{C}$.  But the proof of Lemma 4 shows that $T_2^{-1}(Z)=K^{-1}(Z)$.  It follows that $Z$ lies on $K(\textsf{C})=\textsf{N}$.   This implies that $\Phi_P$ has a fixed point in the ordinary plane, so it must be a homothety.  Now $Z$ is a fixed point of $\Phi_P$ and this map takes $\textsf{N}$ to $\textsf{I}$, so $Z$ also lies on the circle $\textsf{I}$.  Finally, any affine map takes the tangent lines on a conic to the tangent lines on the image conic [2, p. 87], so the tangent $t$ to $\textsf{N}$ at $Z$ gets mapped by $\Phi_P$ to the tangent $t'=\Phi_P(t)$ to $\textsf{I}$ at $Z$. As a homothety, $\Phi_P$ fixes the finite point $Z$ as well as the slope of the line $t$, so $\Phi_P(t) = t$.  In other words, the tangent lines $t$ and $t'$ are identical, which means that  $\textsf{N}$ and $\textsf{I}$ must touch each other at $Z$.  This completes the proof of the Main Theorem.  This proof applies equally to the case $P=P_A$, by the results of Section 3, and yields the following theorem.  \bigskip

\noindent {\bf Theorem 3.} The map $\Phi_A$ of Theorem 2 is a homothety taking $\textsf{N}$ to the excircle $\textsf{E}_A$ and fixing the point $Z_A= (x(y-z)^2,y(x-z)^2,z(x-y)^2)$, where $(x,y,z)$ are homogeneous barycentric coordinates for the external Gergonne point $P_A$ which is opposite vertex $A$.  The circles $\textsf{N}$ and $\textsf{E}_A$ are tangent to each other at $Z_A$.   There are similar statements for the other excircles.  \bigskip

\noindent The Main Theorem and Theorem 3 now imply Feuerbach's Theorem.  A completely synthetic proof of the Main Theorem is given in [11, Section 8], using concepts from projective geometry.

\section{Feuerbach's Theorem in Hilbert Planes.}

\noindent A close examination of the above proof shows that we have only used arguments that are valid in any Hilbert plane $\Pi$ with the parallel postulate.  A Hilbert plane is any set in which points and lines are defined which satisfy Hilbert's basic group of 13 axioms of incidence, betweenness, and congruence.  (See [6], [8].)  In addition, since the parallel postulate holds, $\Pi$ is isomorphic to the cartesian plane $\textsf{F}^2$ over an ordered, Pythagorean field $\textsf{F}$.  An ordered field is one in which a set $\mathcal{P}$ of positive elements is distinguished, such that $\mathcal{P}$ is closed under addition and multiplication; and such that for any element $a \in \textsf{F}$ we have either $a \in \mathcal{P}, a = 0$, or $-a \in \mathcal{P}$.  In such a field the elementary theory of inequalities holds as it does in $\mathbb{R}$, and this allows one to define the notion of betweenness.  The field $\textsf{F}$ is Pythagorean if $a \in \textsf{F}$ implies that  $\sqrt{1+a^2} \in \textsf{F}$.  This condition is required in order to be able to define distance and mark off given distances on given rays, and to prove the SAS criterion for congruence (one of Hilbert's congruence axioms).  The theory of similarity and affine maps, Desargues' Theorem, and Ceva's Theorem are all valid in this generality (see [8, ch. 4] and [6, pp. 462-463]).  The field $\textsf{F}$ may be taken to be the field of equivalence classes of segments in $\Pi$ under the equivalence relation of congruence (together with 0 and negatives of segment classes).  \medskip

\noindent It is also proved in [8, pp. 55, 57-58, 173, 175-176] that the nine-point circle, the incenter, and incircle $\textsf{I}$ are all defined in this situation.  Once we have the incenter, we can drop perpendiculars to the sides to determine the points $D, E,$ and $F$, whose cartesian coordinates are in the field generated by the coordinates of the vertices and the coordinates of the incenter.   This allows us to determine the coordinates of the intersection $P$ of $AD$ and $BE$, for example.  We can then determine $t, u \in \textsf{F}$ for which $v_D=tv_B+(1-t)v_C$ and $v_P=uv_A+(1-u)v_D$.  This requires solving linear equations only, so the barycentric coordinates of $P$ are elements of $\textsf{F}$.  It follows that the coordinates of $Z$ (barycentric and cartesian) lie in $\textsf{F}$, and all of the calculations in sections 4 and 5 are valid.  The same considerations apply to the excircles.  Thus Feuerbach's Theorem holds in the Hilbert plane $\Pi$.  \medskip

\noindent It may seem that in looking at the family of maps $\Phi_P$ in Sections 4 and 5 we are generalizing more than is necessary.  However, for any point $P=(x,y,z)$ with positive barycentric coordinates, there is a triangle for which the Gergonne point of the triangle is $P$ (assuming an additional condition on the field $\textsf{F}$), so that the map $\Phi$ of Section 2 coincides with $\Phi_P$.  As we allow $x, y$, and $z$ to vary over positive elements, we are in essence varying over all possible triangles.  The proof of this claim depends on the following lemma.  \bigskip

\noindent {\bf Lemma 5.} Let $Q $ be the point with homogeneous barycentric coordinates $(x', y', z') = (x(y+z), y(z+x), z(x+y))$. Then: \medskip

a) The intersection of lines $AQ$ and $EF$ is the midpoint of segment $EF$.  Similarly, $BQ \cdot DF$ and $CQ \cdot DE$ are the midpoints of $DF$ and $DE$. (See [7], [13].)  \medskip

b) If $P$ is the Gergonne point of triangle $ABC$, then $Q$ is the incenter, i.e., the center of the incircle of $ABC$. \medskip

c) The incenter $Q$ of $ABC$ has barycentric coordinates $(a, b, c)$, where $a, b$, and $c$ are the lengths of the sides of $ABC$ opposite $A, B, C$, respectively. (See [9], [10].) \medskip

\noindent {\it Proof.}  a) The midpoint of segment $EF$ corresponds to the average of the vectors $v_E$ and $v_F$:

$$\frac{1}{2}v_E+\frac{1}{2}v_F=\frac{1}{2}\bigg(\frac{x}{x+z}v_A+\frac{z}{x+z}v_C\bigg)+\frac{1}{2}\bigg(\frac{x}{x+y}v_A+\frac{y}{x+y}v_B\bigg)$$

$$=\frac{x^2}{(x+y)(x+z)}v_A+\frac{s}{2(x+y)(x+z)}\bigg(\frac{x'}{s}v_A+\frac{y'}{s}v_B+\frac{z'}{s}v_C\bigg) \hspace{.2 in} (s = x'+y'+z')$$

$$=\frac{x^2}{(x+y)(x+z)}v_A+\frac{s}{2(x+y)(x+z)}v_Q.$$ \smallskip

\noindent The last representation implies that the midpoint of $EF$ lies on the line $AQ$ and completes the proof of the claim. \bigskip

\noindent b) By part a), the line $AQ$ passes through the midpoint of $EF$.  Since $AE$ and $AF$ are tangents to $\textsf{I}$, $AQ$ lies on the center of $\textsf{I}$.  The same argument applies to $BQ$ and $CQ$, so $Q$ must be the center of $\textsf{I}$.  \bigskip

\noindent c) We use the fact that in triangle $ABC$ the angle bisector $AQ$ divides the segment $BC$ into segments $BH$ and $HC$ for which $BH/HC=AB/AC=c/b$.  (This is Euclid (VI.3), but see [8, p. 181, Exer. 20.2].)  It follows that the barycentric coordinates of the point $H$ are $(0,b,c)$.  In the same way the intersections of the other angle bisectors with the opposite sides are $(a, 0, c)$ and $(a, b, 0)$.  The claim now follows from (4.6). \bigskip

\noindent To compute the barycentric coordinates of the Gergonne point we solve the system of equations

$$x(y+z)=\lambda a, \quad y(x+z)=\lambda b, \quad z(x+y)=\lambda c$$ \smallskip

\noindent for $x, y, z$.  Distributing on the left sides of these equations and subtracting each equation from the sum of the other two we obtain

$$2yz=\lambda(b+c-a), \quad 2xz = \lambda(a+c-b), \quad 2xy=\lambda(a+b-c).\eqno{(6.1)}$$ \smallskip

\noindent Equation (4.7) shows that the right sides in (6.1) are the barycentric coordinates of the point $P'=(b+c-a,a+c-b,a+b-c)$.  Hence, the Gergonne point $P$ has barycentric coordinates $(x, y, z)$, where

$$x = \frac{1}{b+c-a}, \quad y=\frac{1}{a+c-b}, \quad z=\frac{1}{a+b-c}.\eqno{(6.2)}$$ \smallskip

\noindent For any positive numbers $x, y, z$ there are positive numbers $a, b, c$ for which these equations hold, which is easy to see from the formula for the map $K^{-1}$:

$$K^{-1} {\vecd a b c}={\vecd {1/x} {1/y} {1/z}}. $$ \smallskip

\noindent Thus $(a,b,c) = (1/x, 1/y, 1/z)K$ (since $K$ is a symmetric matrix), and if $x, y, z$ are positive it is immediate that $a, b, c$ are also positive.  This proves the claim made above, since $a, b, c$ certainly form the sides of a suitable triangle.  (This remark requires the underlying field of segments $\textsf{F}$ to be Euclidean, i.e., that it contain the square-roots of its positive elements.  See [6, p. 156, Exer. 5] and [8, p. 144].) \medskip

\noindent The formulas in (6.2) also allow us to compute the barycentric coordinates $(\xi, \eta, \zeta)$ of the Feuerbach point $Z$.  Using (6.2) in the definition of the point $Z$, and multiplying through by the homogeneous factor $\lambda=\frac{1}{4}(b+c-a)^2(a+c-b)^2(a+b-c)^2$, gives

\begin{eqnarray*}
\xi=(b-c)^2(b+c-a),\\
\eta=(a-c)^2(a+c-b),\\
\zeta=(a-b)^2(a+b-c).
\end{eqnarray*}

\medskip

\noindent (Also see [9] and [10].)  Since the cartesian coordinates of the Feuerbach point are the coordinates of the vector

$$v_Z=\frac{\xi}{s}v_A+\frac{\eta}{s}v_B+\frac{\zeta}{s}v_C, \quad s=\xi+\eta+\zeta,$$ \smallskip

\noindent it is clear that the cartesian coordinates of $Z$ lie in the field generated over $\mathbb{Q}$ by the cartesian coordinates of the vertices $A, B, C$ and the coordinates of the incenter $Q$.  In particular, if the vertices have coordinates in the Hilbert field $\Omega$, the smallest ordered Pythagorean field (see [8, p. 145]), the coordinates of $Z$ lie in $\Omega$ as well.  \medskip

\noindent One can also show that if the coordinates of the vertices lie in $\mathbb{Q}$, then the coordinates of the points $Q$ and $Z$ lie in a quartic extension of $\mathbb{Q}$.  This is because the equation of the line $AB$ can be written as

$$(b_2-a_2)x - (b_1-a_1)y + r=0,$$ \smallskip

\noindent so the distance of $Q=(x,y)$ from $AB$ has the form

$$d=\frac{|(b_2-a_2)x - (b_1-a_1)y + r|}{\sqrt{(b_2-a_2)^2+(b_1-a_1)^2}}=\frac{|(b_2-a_2)x - (b_1-a_1)y + r|}{c},$$ \smallskip

\noindent since the length of $AB$ is $c$.  Setting this expression equal to the distance of $Q$ from the line $AC$, squaring, and factoring the resulting expression, shows that the equation of the angle bisector of $\angle B$ has coordinates in the field $\mathbb{Q}(a/c)$.  Doing the same for angle $\angle C$ gives that $Q$ has coordinates in the field $\mathbb{Q}(a/c, a/b)$, which is generally a quartic extension of $\mathbb{Q}$. 

\bigskip

\noindent {\bf Example.} Consider triangle $ABC$, where

$$v_A=(4,0), \quad v_B=(0, 3), \quad v_C=(-2, 0).$$

\noindent The lengths of the sides are

$$a=\sqrt{13}, \quad b=6, \quad c = 5.$$

\noindent It is not hard to verify that the equation of the nine-point circle is

$$\textsf{N}: \quad x^2-x+y^2-\frac{17}{6}y=0.$$ \smallskip

\noindent To find the incenter we find the equations of the angle bisectors of $\angle A$ and $\angle C$ using the formula for the distance of a point from a line.  We find:

$$\textrm{bisector of} \hspace{.05 in} \angle A: x+3y-4=0, \quad \textrm{bisector of} \hspace{.05 in} \angle C:  3x-(2+\sqrt{13})y+6=0.$$ \smallskip

\noindent Thus the incenter is given by

$$v_Q=\left( \frac{-3+\sqrt{13}}{2}, \frac{11-\sqrt{13}}{6} \right)$$

\noindent and the equation of the incircle is

$$\textsf{I}: \quad x^2+(3-\sqrt{13})x+y^2+\left(\frac{-11+\sqrt{13}}{3} \right)y+\frac{11-3\sqrt{13}}{2}=0.$$

\noindent Solving the equations for $\textsf{N}$ and $\textsf{I}$ simultaneously gives the cartesian coordinates

$$v_Z=\left( \frac{-79+5\sqrt{13}}{102}, \frac{56-10\sqrt{13}}{51} \right) \eqno{(6.3)}$$ \smallskip

\noindent of the Feuerbach point $Z$ of $ABC$.  The coordinates of the points of tangency of the incircle with the sides and their corresponding cevians are

$$v_D=\left( \frac{-13+\sqrt{13}}{13}, \frac{39+3\sqrt{13}}{26} \right), \quad AD: y=-\left( \frac{11+\sqrt{13}}{36} \right)(x-4);$$

$$v_E=\left( \frac{-3+\sqrt{13}}{2},0 \right), \quad BE: y=-\frac{3}{2}(3+\sqrt{13})x+3;$$

$$v_F=\left( \frac{-2+2\sqrt{13}}{5}, \frac{33-3\sqrt{13}}{10} \right), \quad CF: y= \frac{57-15\sqrt{13}}{4}(x+2).$$ \smallskip

\noindent Hence the Gergonne point $P$ is given by

$$v_P=\left( \frac{-115+37\sqrt{13}}{127}, \frac{177+6\sqrt{13}}{127} \right).\eqno{(6.4)}$$ \smallskip

\noindent Using (6.2) we can find barycentric coordinates for $P$ directly:

$$P \rightarrow \left( \frac{11+\sqrt{13}}{108}, \frac{1+\sqrt{13}}{12}, \frac{-1+\sqrt{13}}{12} \right),$$

\noindent and verify that these coordinates yield the same result (6.4) for the cartesian coordinates of $P$.  Similarly, the formulas for $\xi, \eta, \zeta$ yield

$$\xi = 11-\sqrt{13}, \quad \eta =-168+48\sqrt{13} , \quad \zeta = -107+37\sqrt{13},$$

$$s=\xi + \eta + \zeta =-264+ 84\sqrt{13},$$ \smallskip

\noindent for the barycentric coordinates of $Z$, and we check directly that

$$\frac{\xi}{s}v_A+\frac{\eta}{s}v_B+\frac{\zeta}{s}v_C=$$

$$\left( \frac{151+55\sqrt{13}}{1836} \right) v_A+\left( \frac{56-10\sqrt{13}}{153} \right) v_B+\left( \frac{1013+65\sqrt{13}}{1836} \right) v_C= v_Z$$

\noindent coincides with the result of (6.3).

\section{References.}

\noindent [1] N. Altshiller-Court, {\it College Geometry, An Introduction to the Modern Geometry of the Triangle and the Circle}, reprinted by Dover Publications, 2008.  \medskip

\noindent [2] D.A. Brannan, M.F. Esplen, and J.J. Gray, {\it Geometry}, Cambridge University Press, 1999. \medskip

\noindent [3] H.S.M. Coxeter, {\it The Real Projective Plane}, 3rd edition, Springer, 1993. \medskip

\noindent [4] H.S.M. Coxeter, {\it Introduction to Geometry}, 2nd edition, Wiley, 1969. \medskip

\noindent [5] H. Eves, {\it College Geometry}, Jones and Bartlett Publishers, Boston, 1995. \medskip

\noindent [6] M. Greenberg, {\it Euclidean and Non-Euclidean Geometries: Development and History}, 4th edition, W. H. Freeman, 2008. \medskip

\noindent [7] D. Grinberg, Hyacinthos Message 6423 (2003), http://tech.groups.yahoo.com/ \smallskip

\noindent group/Hyacinthos.  \medskip

\noindent [8] R. Hartshorne, {\it Geometry: Euclid and Beyond}, Springer, 2000. \medskip

\noindent [9] C. Kimberling, Central points and central lines in the plane of a triangle, Mathematics Magazine 67 (1994), 163-187.  \medskip

\noindent [10] C. Kimberling, Encyclopedia of Triangle Centers, http://faculty.evansville.edu \smallskip

\noindent /ck6/encyclopedia/. \medskip

\noindent [11] I. Minevich and P. Morton, Synthetic cevian geometry, IUPUI preprint 2009, available as PR 09-01 at http://www.math.iupui.edu/research/preprint/
2009/pr09-01.pdf. \medskip

\noindent [12] P. Yiu, A Tour of Triangle Geometry, at www.math.fau.edu/yiu/
TourOfTriangleGeometry /MAAFlorida37040428.pdf. \medskip

\noindent [13] P. Yiu, Hyacinthos Message 1790 (2000), http://tech.groups.yahoo.com/ \smallskip

\noindent group/Hyacinthos. \medskip

\bigskip

\noindent Dept.of Mathematical Sciences, LD 270 \smallskip

\noindent Indiana University Ð Purdue University at Indianapolis (IUPUI) \smallskip

\noindent 402 N. Blackford St. \smallskip

\noindent Indianapolis, IN 46202 \smallskip

\noindent {\it e-mail}: pmorton@math.iupui.edu

\end{document}